\documentclass{ifacconf}

\renewcommand\tilde{\widetilde}

\newtheorem{nnassumption}{\bf Assumption}

\newtheorem{nntheorem}{\bf Theorem}
\newenvironment{theorem}{\begin{nntheorem}\it}{\end{nntheorem}}
\newtheorem{nndefinition}{\bf Definition}
\newenvironment{definition}{\begin{nndefinition}\it}{\end{nndefinition}}
\newtheorem{nnproposition}{\bf Proposition}
\newenvironment{proposition}{\begin{nnproposition}\it}{\end{nnproposition}}
\newtheorem{nnproblem}{\bf Problem}

\newtheorem{nnlemma}{\bf Lemma}

\newtheorem{nnremark}{\bf Remark}
\newenvironment{remark}{\begin{nnremark}}{\hfill \hspace*{1pt}\hfill $\circ$\end{nnremark}}
\newenvironment{proof}{{\bf Proof.}}{\hfill \hspace*{1pt}\hfill $\bullet$}

\newtheorem{example}{Example}

\usepackage{amsfonts,latexsym,amsmath,amssymb}
\usepackage[mathscr]{eucal}
\usepackage{graphicx,color}
\usepackage{comment}

\usepackage{graphicx}      
\usepackage{natbib}        
\usepackage[utf8]{inputenc}
\begin{document}
\begin{frontmatter}

\title{{Stability analysis of a 1D wave equation with a nonmonotone distributed damping\thanksref{footnoteinfo}}}

\thanks[footnoteinfo]{This works has been partially supported by Advanced Grant DYCON (Dynamic Control) of the European Research Council Executive Agency.}

\author[First,Four,Five]{Swann Marx}, 
\author[Second]{Yacine Chitour}, 
\author[Third]{Christophe Prieur}

\address[First]{LAAS-CNRS, Universit\'e de Toulouse, CNRS, 7 avenue du colonel Roche, 31400, Toulouse, France
        {\tt\small marx.swann@gmail.com}.}
\address[Second]{Laboratoire des Signaux et Syst\`emes (L2S), CNRS - CentraleSupelec - Universit\'e Paris-Sud, 3, rue Joliot Curie, 91192, Gif-sur-Yvette, France,
        {\tt\small yacine.chitour@lss.supelec.fr}.}
\address[Third]{Univ. Grenoble Alpes, CNRS, Grenoble INP, Gipsa-lab, F-38000 Grenoble, France,
        {\tt\small christophe.prieur@gipsa-lab.fr}}
\address[Four]{DeustoTech, University of Deusto, 48007 Bilbao, Basque Country,
Spain.}
\address[Five]{Facultad de Ingenieria, Universidad de Deusto, Avda.Universidades,
24, 48007, Bilbao - Basque Country – Spain.}

\begin{abstract}                
This paper is concerned with the asymptotic stability analysis of a one dimensional wave equation subject to a nonmonotone distributed damping. A well-posedness result is provided together with a precise characterization of the asymptotic behavior of the trajectories of the system under consideration. The well-posedness is proved in the nonstandard $L^p$ functional spaces, with $p\in [2,\infty]$, and relies mostly on some results collected in \cite{haraux1D}. The asymptotic behavior analysis is based on an attractivity result on a specific infinite-dimensional linear time-variant system. 
\end{abstract}

\begin{keyword}
1D wave equation, nonlinear control, Lyapunov functionals.
\end{keyword}

\end{frontmatter}

\section{Introduction}


This paper is concerned with the asymptotic behavior of a one-dimensional wave equation subject to a nonmonotone nonlinear damping. For control systems, it is always crucial to consider nonlinear feedback laws. In particular, the definition of the nonlinearity under consideration in this paper includes the saturation, which models amplitude limitations on the actuator. Such a phenomenon appears in most of control systems, and it can lead to undesirable behavior in terms of stability. Moreover, considering that such a nonlinearity is nonmonotone is crucial for control systems, because the model of the nonlinearity might have some error.

However, as illustrated in many papers (\cite{alabau2012some}, \cite{haraux1D}, \cite{marx2018stability}, \cite{slemrod1989mcss}, \cite{feireisl1993strong}, etc.), the monotone property of the nonlinearity is crucial to first prove the asymptotic stability of the system under consideration and second characterize the trajectory of the latter system. In order to characterize the asymptotic behavior of the trajectory of the system under consideration, we consider the initial conditions in an another functional setting than the classical one, that is $L^p(0,1)$, with $p\in [2,\infty]$.

There exists a vast litterature about linear PDEs subject to monotone nonlinear dampings. For instance, in \cite{slemrod1989mcss}, asymptotic stability of the origin of abstract control systems subject to monotone nonlinear dampings is proved, using an infinite-dimensional version of LaSalle's invariance principle. These results have been then extended to more general infinite-dimensional systems in \cite{map2017mcss}. More recently, in \cite{marx2018stability}, trajectory of such systems have been characterized via Lyapunov techniques. Let us also mention \cite{prieur2016wavecone} and \cite{mcpa2017siam}, where a wave equation and a nonlinear Korteweg-de Vries, respectively, subject to a nonlinear monotone damping are considered and where the global asymptotic stability is proved. Finally, in \cite{kang2017boundary}, the local asymptotic stability of a heat equation coupled with an ODE, controlled from the boundary with a saturated feedback law, is proved.

The case of nonmonotone nonlinear damping have been considered in few papers. In \cite{feireisl1993strong}, the global asymptotic stability of a one-dimensional wave equation subject to a nonmonotone nonlinear damping is proved, thanks to some compensated compactness technique. The characterization of the trajectories is, however, not provided. In \cite{MVC}, the trajectory of a wave equation in two dimensions suject to a nonmonotone damping is characterized, but for only a specific nonmonotone damping. In this paper, we rather focus on a more general nonlinear nonmonotone damping, but only for a one dimensional wave equation. We are able to characterize the trajectories of the system by studying the functional spaces $L^p(0,1)$.

There exists also few papers dealing with the one dimensional wave equation in this functional setting. Let us mention \cite{haraux1D}, which derives a well-posedness analysis of a one dimensional wave equation in this functional setting. Moreover, an optimal decay rate is obtained for this equation. Recently, in \cite{amadori2019decay}, a similar result have been obtained, using techniques coming from conservation laws theory. Note that both of these results hold true only for monotone nonlinear damping, which is not the case of our paper.

In our paper, after introducing a general nonlinear nonmonotone damping, we propose a well-posedness analysis of a one-dimensional wave equation subject to such a damping. Futhermore, using Lyapunov techniques, we are able to characterize the trajectory of this system in the $L^p$, $p\in [2,\infty)$, provided that the initial conditions are in $L^\infty$. This proof is mainly based on a result about linear time-varying infinite-dimensional system, which is introduced and proved in Appendix \ref{sec_appendix}.

This paper is organized as follows. In Section \ref{sec_main}, the main results of the paper are collected. We propose a well-posedness and an asymptotic stability theorems. The second result proposes futhermore a precise characterization of the trajectories of the system. Section \ref{sec_proof} is devoted to the proof of the main results. Finally, Section~\ref{sec_conclusion} collects some concluding remarks, together with further research lines to be followed. Appendix \ref{sec_appendix} introduces a result of independant interest about exponential convergence of a specific time-varying infinite-dimensional linear systems, following a finite-dimensional strategy provided in \cite{chitour1995continuity} and \cite{liu1996finite}, but which is instrumental for the proof of our asymptotic stability theorems.

\textbf{Notation:} For any $p\in [2,\infty)$, the space $L^p(0,1)$ denotes the space of functions $f$ satisfying $\left(\int_0^1 |f(x)|^p dx\right)^{\frac{1}{p}}<+\infty$. The space $L^\infty(0,1)$ denotes the space of functions satisfying $\mathrm{ess}\sup_{x\in [0,1]} |f(x)|\leq +\infty$. For any $p\in [2,\infty]$, the Sobolev space $W^{1,p}(0,1)$ (resp. $W^{2,p}(0,1)$) is defined as follows $W_0^{1,p}(0,1):=\lbrace f\in L^p(0,1)\mid f^\prime\in L^p(0,1)\text{ and }f(0)=f(1)=0\rbrace$ (resp. $W^{2,p}(0,1):=\lbrace f\in L^p(0,1)\mid f^\prime,f^{\prime\prime}\in L^p(0,1)\rbrace$).   


\section{Main results}
\label{sec_main}

The aim of this paper is to provide an asymptotic stability analysis of the following system:
\begin{equation}
\label{wave-nonmonotone}
\left\{
\begin{split}
&z_{tt}(t,x) = z_{xx}(t,x)-\sigma(a(x)z_t(t,x)),\: (t,x)\in \mathbb{R}_+\times [0,1]\\
&z(t,0)=z(t,1) = 0,\: t\in\mathbb{R}_+\\
&z(0,x) = z_0(x),\: z_t(0,x)=z_1(x),\: x\in [0,1],
\end{split}
\right.
\end{equation}
where $z$ denotes the state, $a:[0,1]\to \mathbb{R}_+$ is measurable and bounded by some positive constant $a_\infty$ and $\sigma$ is a scalar nonlinear damping which satisfies the following properties:
\begin{definition}[Scalar nonlinear damping]
\label{def-damping}
A function $\sigma:\mathbb{R}\rightarrow \mathbb{R}$ is said to be a \emph{scalar damping function} if 
\begin{itemize}
\item[1.] It is locally Lipschitz and odd;
\item[2.] One has $\sigma(0)=0$;
\item[3.] For any $s\in \mathbb{R}$, $\sigma(s)s>0$;
\item[4.] The function $\sigma$ is differentiable at $s=0$ with $\sigma'(0)=C_1$ for some $C_1>0$.
\end{itemize}
\end{definition}
Due to this definition (especially, item 2), the origin is an equilibrium point for \eqref{wave-nonmonotone}. Note that this nonlinearity is not assumed to be \emph{monotone}. The monotone property is in many cases really useful for either the well-posedness of the equation or the asymptotic stability of the origin. We refer the reader to \cite{marx2018stability} for a complete discussion on this topic.

\begin{figure}[h!]
\includegraphics[scale=0.5]{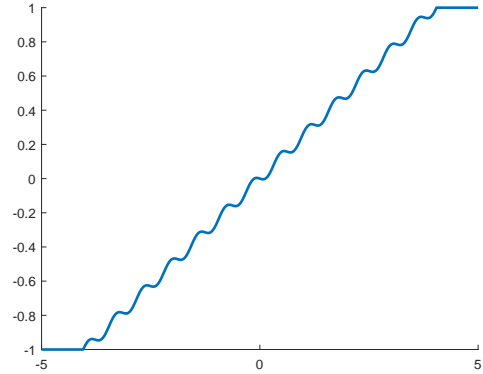}
\caption{For any $s\in [-5,5]$, the figure illustrates the function $\sigma$ given by \eqref{example-nonmonotone}}
\label{fig:nonmonotone}
\end{figure}

\begin{example}[Example of nonlinear dampings]
Below are listed some examples of nonlinear dampings:
\begin{itemize}
\item[1.] The classical saturation, defined as follows:
\begin{equation}
\sigma(s)=\mathrm{sat}(s):=\left\{
\begin{split}
&s  &\text{ if } |s|\leq 1,\\
&\frac{s}{|s|}  &\text{ if } |s|\geq 1,
\end{split}
\right.
\end{equation}
satisfies all the properties of Definition \ref{def-damping}.
\item[2.] The following nonlinearity
\begin{equation}
\label{example-nonmonotone}
\sigma(s) = \mathrm{sat}\left(\frac{1}{4}s-\frac{1}{30}\sin(10s)\right)
\end{equation}
is also a nonlinear damping. Note moreover that it is not monotone, as illustrated by Figure \ref{fig:nonmonotone}.
\end{itemize}
\end{example}

As illustrated in \cite{marx2018stability}, some regularity is needed to obtain a characterization of the asymptotic stability of \eqref{wave-nonmonotone}. To be more precise, we need the state $z_t$ to be bounded in $L^\infty(0,1)$. With a monotone nonlinearity $\sigma$, one would have this regularity result thanks to some nonlinear semigroup theorems. In the case of the system under consideration in this paper, we need to follow another strategy. 

Our strategy relies on the introduction of these functional spaces:
\begin{equation}
\begin{split}
&H_p(0,1):= W_0^{1,p}(0,1)\times L^p(0,1)\\
&D_p(0,1):=W^{2,p}(0,1)\cap W_0^{1,p}(0,1)\times W_0^{1,p}(0,1),
\end{split}
\end{equation}
with $p\in [2,\infty]$. The first functional space is equipped with the following norm:
\begin{equation}
\begin{split}
\Vert (z,z_t)\Vert_{H_p(0,1)} := &\left(\int_0^1 [z_x(x)|^p dx\right)^{\frac{1}{p}} \\
&+ \left(\int_0^1 |z_t(x)|^p dx\right)^{\frac{1}{p}},\: \forall p\in [2,\infty),\\
\Vert (z,z_t)\Vert_{H_\infty(0,1)}:=&\Vert z_x\Vert_{L^\infty(0,1)} + \Vert z_t\Vert_{L^\infty(0,1)},\: \text{for $p=\infty$.}
\end{split}
\end{equation}
The second one is equipped with the following norm
\begin{equation}
\begin{split}
\Vert (z,z_t)\Vert_{D_p(0,1)}:=&\left(\int_0^1 |z_{xx}(x)|^p dx\right)^{\frac{1}{p}} \\
&+ \left(\int_0^1 |z_{tx}(x)|^p dx\right)^{\frac{1}{p}},\: \forall p\in [2,\infty)\\
\Vert (z,z_t)\Vert_{D_\infty(0,1)}:=&\Vert z_{xx}\Vert_{L^\infty(0,1)} + \Vert z_{t,x}\Vert_{L^\infty(0,1)},\: \text{for $p=\infty$.}
\end{split}
\end{equation}
Considering that $(z_0,z_1)\in H_\infty(0,1)$, our aim is to prove an $L^\infty$ regularity on the state $z_t$ so that we can obtain a characterization of the asymptotic behavior of the trajectory in $H_p(0,1)$, with $p\in [2,\infty)$. But before stating such a result, we need a suitable notion of solution for systems modeled with \eqref{wave-nonmonotone}. Such a notion is provided by the following well-posedness theorem:
\begin{theorem}[Well-posedness]
\label{thm-wp}
For any initial conditions $(z_0,z_1)\in H_\infty(0,1)$, there exists a unique solution $(z,z_t)\in L^\infty(\mathbb{R}_+;W^{1,\infty}(0,1))\times W^{1,\infty}(\mathbb{R}_+;L^\infty(0,1))$ to \eqref{wave-nonmonotone}. Moreover, the following inequality is satisfied, for all $t\geq 0$
\begin{equation}
\Vert (z,z_t)\Vert_{H_\infty(0,1)}\leq 2\max(\Vert z_0^\prime\Vert_{L^\infty(0,1)}, \Vert z_1\Vert_{L^\infty(0,1)}).\footnote{The term $z_0^\prime$ denotes the derivative of $z_0$ with respect to $x$.}
\end{equation}
\end{theorem}
Now that the functional setting is introduced, we are in position to state our asymptotic stability result
\begin{theorem}[Semi-global exponential stability]
\label{thm-as}
Consider initial conditions $(z_0,z_1)\in H_\infty(0,1)$ satisfying:
\begin{equation}
\Vert (z_0,z_1)\Vert_{H_\infty(0,1)}\leq R,
\end{equation}
where $R$ is a positive constant. Then, for any $p\in [2,\infty)$, there exist two positive constants $K:=K(R)$ and $\beta:=\beta(R)$ such that
\begin{equation}
\Vert (z,z_t)\Vert_{H_p(0,1)}\leq Ke^{-\beta t}\Vert (z_0,z_1)\Vert_{H_p(0,1)},\quad \forall t\geq 0.
\end{equation}
\end{theorem}

\section{Proof of the main results}

\label{sec_proof}

\subsection{Proof of Theorem \ref{thm-wp}}

The proof of Theorem \ref{thm-wp} relies on the following wave equation with a source term
\begin{equation}
\label{wave-source}
\left\{
\begin{split}
&z_{tt}(t,x) = z_{xx}(t,x) +h(t,x),\: (t,x)\in\mathbb{R}_+\times [0,1]\\
&z(t,0)=z(t,1) = 0,\: t\in\mathbb{R}_+\\
&z(0,x)= z_0(x),\: z_t(0,x) = z_1(x),\: x\in [0,1],
\end{split}
\right.
\end{equation}
where $h$ denotes the source term. From \cite[Theorem 1.3.8]{HarauxBook}, we know that, provided that $z_0,z_1\in H_2$ and that $h\in L^2(\mathbb{R}_+;L^2(0,1))$, there exists a unique solution $z\in C(\mathbb{R}_+;H_0^1(0,1))\cap C^1(\mathbb{R}_+;L^2(0,1))$ to \eqref{wave-source}. In particular, since $H_\infty(0,1)\subset H_2$, this result holds true also for initial conditions $(z_0,z_1)\in H_{\infty}$. The first step of our analysis in this section is to prove that picking initial conditions in $H_\infty(0,1)$ and the source term $h\in L^2(\mathbb{R}_+;L^\infty(0,1))$ improves also the regularity of the solution $z$ itself. 

To do so, our aim is to give an explicit formula for the latter equation, using the reflection method surveyed in \cite{strauss1992partial}. Roughly speaking, this method consists in extending the explicit formulation of trajectory of the wave equation in an unbounded domain to a bounded domain. To do so, we extend the initial datas to the whole line to be odd with respect to both $x=0$ and $x=1$, that is
\begin{equation}
\tilde{z}_0(-x) = -\tilde{z}_0(x)\text{ and } \tilde{z}_0(2-x)=-\tilde{z}_0(x),
\end{equation}
where $\tilde{z}_0$ denotes the $2$-periodic odd extension of $z_0$. A way to do this is to define $\tilde{z}_0$ as follows: 
\begin{equation}
\tilde{z}_0(x) = \left\{
\begin{split}
&z_0(x),\: 0<x<1\\
&-z_0(-x),\: -l<x<0\\
&\text{extended to be of period $2$} 
\end{split}
\right.
\end{equation}
We can define similarly a $2$-periodic odd extension of $z_1$ (resp. $h$), denoted by $\tilde{z}_1$ (resp. $\tilde{h}$). Thanks to \cite[Theorem 1, Page 69]{strauss1992partial}, we can therefore define the explicit trajectory $z$ of \eqref{wave-source} (known as the D'Alembert formula) as follows:
\begin{equation}
\begin{split}
z(t,x) = &\frac{1}{2} \left[\tilde{z}_0(x+t)+\tilde{z}_0(x-t)\right] + \frac{1}{2}\int_{x-ct}^{x+ct}\tilde{z}_1(s) ds \\
&+ \frac{1}{2}\int_0^t\int_{x-(t-s)}^{x+(t-s)} \tilde{h}(w,s)dwds. 
\end{split}
\end{equation}
We can further define $z_t$ as follows
\begin{equation}
\begin{split}
z_t(t,x) =& \frac{1}{2} (\tilde{z}_0^\prime(x+t)-\tilde{z}_0^\prime(x-t)) \\
&+ \frac{1}{2}\left(\tilde{z}_1(x+t) - \tilde{z}_1(x-t)\right)\\
&+\frac{1}{2}\int_0^t \left(\tilde{h}(s,x+(t-s)) -\tilde{h}(s,x-(t-s))\right) ds 
\end{split}
\end{equation}
It is clear from these two latter equations that, when picking $(z_0,z_1)\in H_\infty(0,1)$ and $h\in L^2(\mathbb{R}_+;L^\infty(0,1))$, then $$z\in C(\mathbb{R}_+;W^{1,\infty}(0,1))\cap C^1(\mathbb{R}_+;L^\infty(0,1)).$$
We assume now that $\tilde{h}$ is written as follows
\begin{equation}
\label{source-term}
\tilde{h}(t,x):=-\sigma(\tilde{a}(x)y(t,x)),
\end{equation}
where $\tilde{a}$ is a $2$-periodic extension of $a$, $y\in L^2(0,T;L^2(0,1))$ and $\sigma$ is a scalar nonlinear damping (which is odd, due to Item 1 of Definition \ref{def-damping}). In particular, it means that $h\in L^2(\mathbb{R}_+;L^2(0,1))$. The proof of Theorem \ref{thm-wp} consists first in applying a fixed-point theorem, which will allow us to prove the well-posedness of \eqref{wave-nonmonotone} for a small time $T>0$ and second in using a stability result in \cite{haraux1D}, stated as follows
\begin{proposition}
Let us consider initial condition $z_0,z_1\in H_\infty(0,1)$. If there exists a solution to \eqref{wave-nonmonotone}, then the time derivative of the following functional along the trajectories of \eqref{wave-nonmonotone}
\begin{equation}
\phi(z,z_t):=\int_{0}^1 \left(F(z_t+z_x)+F(z_t-z_x)\right)dx,
\end{equation}
with $F$ any even and convex function, satisfies
\begin{equation*}
\frac{d}{dt}\phi(z,z_t)\leq 0.
\end{equation*}
\end{proposition}
The latter proposition implies that
\begin{equation}
\phi(z,z_t)\leq \phi(z_0,z_1),\quad \forall t\geq 0.
\end{equation}
In particular, following the discussion in the proof of Corollary 2.3. in \cite{haraux1D}, if one picks \footnote{The function $\mathrm{Pos}\: : s\in\mathbb{R}\rightarrow \mathrm{Pos}(s)\in\mathbb{R}_+$ is defined as follows:
\begin{equation*}
\mathrm{Pos}(s)=\left\{
\begin{split}
&s\text{ if } s>0\\
&0\text{ if } s\leq 0.
\end{split}
\right.
\end{equation*}
} $F(s)=[\mathrm{Pos}(|s|-2\max(\Vert z_0\Vert_{L^\infty(0,1)}, \Vert z_1\Vert_{L^\infty(0,1)}))]^2$, then one obtains that
\begin{equation}
\phi(z,z_t)=0,
\end{equation}
which implies that, for all $t\geq 0$
\begin{equation}
\begin{split}
\max(\Vert z_x(t,\cdot)\Vert_{L^\infty(0,1)},&\Vert z_t(t,\cdot)\Vert_{L^\infty(0,1)})\leq \\
&2\max(\Vert z_0^\prime\Vert_{L^\infty(0,1)}, \Vert z_1\Vert_{L^\infty(0,1)}).
\end{split}
\end{equation}
Noticing that $$\Vert (z,z_t)\Vert_{H_\infty(0,1)}\leq \max(\Vert z_x(t,\cdot)\Vert_{L^\infty(0,1)},\Vert z_t(t,\cdot)\Vert_{L^\infty(0,1)}),$$ it is clear then that, for all $t\geq 0$
\begin{equation}
\label{estimate-Haraux}
\Vert (z,z_t)\Vert_{H_\infty(0,1)} \leq 2\max(\Vert z_0^\prime \Vert_{L^\infty(0,1)}, \Vert z_1\Vert_{L^\infty(0,1)})
\end{equation}
This estimate implies that, once one is able to prove that there exists a solution $(z,z_t)$ of \eqref{wave-nonmonotone} in $$L^\infty([0,T];W^{1,\infty}(0,1))\cap W^{1,\infty}([0,T];L^\infty(0,1)),$$ for a small time $T>0$, then the well-posedness of \eqref{wave-nonmonotone} is ensured in $L^\infty(\mathbb{R}_+;W^{1,\infty}(0,1))\cap W^{1,\infty}(\mathbb{R}_+;L^\infty(0,1))$.

Now, we are in position to prove Theorem \ref{thm-wp}.

\begin{proof}
Let us define $\mathcal{F}_T$ the space of measurable functions defined on $[0,T]\times\mathbb{R}$ which are bounded, odd and $2$-periodic in space. We endow $\mathcal{F}_T$ with the $L^\infty$-norm so that it becomes a Banach space. Hence, denoting by $\Vert \cdot\Vert_{T}$ the norm of the latter functional space, we have, for every $y\in\mathcal{F}_T$
\begin{equation}
\Vert y\Vert_T:=\sup_{(t,x)\in [0,T]\times\mathbb{R}} |y(t,x)|.
\end{equation}
Let us consider $\mathbf{B}_K(y)$ the closed ball in $\mathcal{F}_T$ centered at $y\in\mathcal{F}_T$ of radius $K\geq 0$, where $K$ remains to be defined. 

We can now define the mapping with which we will apply a fixed-point
\begin{equation}
\begin{split}
\phi_T : \mathcal{F}_T &\rightarrow \mathcal{F}_T\\
y &\mapsto \phi_T(y),
\end{split}
\end{equation}
where 
\begin{equation}
\begin{split}
&\phi_T(y) = \frac{1}{2}\left(\tilde{z}^\prime_0(x+t) - \tilde{z}_0^\prime(x-t)\right) + \frac{1}{2}\left(\tilde{z}_1(x+t)-\tilde{z}_1(x-t)\right)\\
& -\frac{1}{2}\int_0^t \left(\sqrt{\tilde{a}(x+t-s)}\sigma(\sqrt{\tilde{a}(x+t-s)}y(s,x-(t-s))\right.\\ 
&\left.+ \sqrt{\tilde{a}(-x+t-s)}\sigma(\sqrt{\tilde{a}(-x+t-s)}y(s,-x-(t-s))\right)ds
\end{split}
\end{equation}
Using the fact that $\sigma$ is locally Lipschitz, note that, for every $y_0\in\mathcal{F}_T$ and $K>0$, there exists a positive constant $C(y_0,K)$ such that, for every $y,\tilde{y}\in \mathbf{B}_K(y_0)$ 
\begin{equation}
\label{contraction}
\Vert \phi_T(y)-\phi_T(\tilde{y})\Vert_T\leq C(y_0,K)T\Vert y-\tilde{y}\Vert_{T}.
\end{equation}
Pick $K$ such that
\begin{equation}
K:=2\left(\Vert \phi_T(y_0)\Vert_T + \Vert y_0\Vert_T\right)
\end{equation}
and $T$ sufficiently small such that
\begin{equation}
\label{choice-T}
C(y_0,K)T\leq 1.
\end{equation}
Consider a sequence $(y_n)_{n\in\mathbb{N}}$ defined as 
\begin{equation*}
y_{n+1} = \phi_T(y_n).
\end{equation*}

If we prove that this sequence converges to $y^\star$, then we can deduce that
\begin{equation*}
y^\star = \phi_T(y^\star).
\end{equation*}
This means that there exists a fixed-point for the mapping $\phi_T$, which implies in particular that \eqref{wave-nonmonotone} is well-posed in the desired functional spaces. To prove the convergence of this sequence, we prove that it is a Cauchy sequence. Indeed, since $\mathbf{B}_K(y_0)$ is complete, any Cauchy sequence converges in this set. By induction, one can prove that
\begin{equation}
\label{inequality-contraction}
\Vert y_{n+1}-y_{n}\Vert_T\leq (C(y_0,K)T)^n \Vert y_0-y_1\Vert_T
\end{equation} 
and
\begin{equation}
\label{ball-property}
y_n\in \mathbf{B}_K(y_0).
\end{equation}
Indeed, thanks to the choice of $T$ in \eqref{choice-T}, these two properties are easily proved for $n=0$ and, moreover, we have with \eqref{contraction}, for all $n\geq 1$
\begin{equation}
\begin{split}
\Vert y_{n+1}-y_n\Vert_T=&\Vert \phi(y_n)-\phi(y_{n-1})\Vert_T\\
\leq &C(y_0,K)T\Vert \phi(y_{n-1})-\phi(y_{n-2})\Vert_T\\
\leq &(C(y_0,K)T)^{n}\Vert y_0-y_1\Vert_T.
\end{split} 
\end{equation}
The inequality \eqref{inequality-contraction} can be deduced from the above inequality. The property \eqref{ball-property} can be proved as follows:
\begin{equation}
\begin{split}
\Vert y_{n+1}\Vert \leq & \Vert y_n\Vert_T + \Vert y_1\Vert_T + \Vert y_0\Vert_T\\
\leq & \Vert y_0\Vert_T + \Vert \phi_T(y_0)\Vert_T + \Vert \phi_T(y_0)\Vert_T + \Vert y_0\Vert_T\\
\leq & K
\end{split}
\end{equation}
where in the first line we have used \eqref{inequality-contraction} and, in the second line, we have used the fact that $y_n\in\mathbf{B}_K(y_0)$. 

The two properties \eqref{inequality-contraction} and \eqref{ball-property} show that the sequence $y_n$ is a Cauchy sequence. Since $\mathbf{B}_K(y_0)$ is a complete set, this sequence is therefore convergent. This means in particular that there exists a fixed-point to the mapping $\phi_T(y_0)$, which implies that, for sufficiently small time $T$, there exists a unique solution $z\in L^\infty(0,T;W^{1,\infty}(0,1))$ and $z_t\in L^\infty(0,T;L^\infty(0,1))$. Thanks to \eqref{estimate-Haraux}, we can deduce that there exists unique solution $z\in L^\infty(\mathbb{R}_+;W^{1,\infty}(0,1))$ and $z_t\in L^\infty(\mathbb{R}_+;L^\infty(0,1))$. This concludes the proof of Theorem \ref{thm-wp}.
\end{proof}

\subsection{Proof of Theorem \ref{thm-as}}
\label{sec_as}

The proof of Theorem \ref{thm-as} is divided into two steps: first, we transform the system \eqref{wave-nonmonotone} into a system in the form \eqref{ltv-system}, which is in particular a linear time-variant system, and apply Theorem \ref{thm-ltv} for the space $H_2$, which is indeed the only Hilbert space among all the spaces $H_p(0,1)$. Second, using the fact that the solutions are bounded in $H_\infty(0,1)$ thanks to Theorem \ref{thm-wp}, and invoking an interpolation theorem, which is the Riesz-Thorin theorem, we conclude. 

\begin{proof}

\textbf{First step: semi-global exponential stability in $H_2$.}
We first fix $p=2$, but still consider the initial conditions $(z_0,z_1)\in H_\infty(0,1)$. Moreover, we suppose that there exists a positive constant $R$ such that the initial conditions satisfy
\begin{equation}
\Vert (z_0,z_1)\Vert_{H_\infty(0,1)} \leq R.
\end{equation}
In particular, due to Theorem \ref{thm-wp}, one has
\begin{equation}
\Vert (z,z_t)\Vert_{H_\infty(0,1)} \leq 2R,\quad \forall t\geq 0.
\end{equation}
The system \eqref{wave-nonmonotone} may be rewritten as follows
\begin{equation}
\label{wave-ltv}
\left\{
\begin{split}
&z_{tt} = z_{xx}-d(t,x)a(x)z_t,\\
&z(t,0) = z(t,1) = 0,\\
&z(0,x) = z_0(x),\: z_t(0,x) = z_1(x),
\end{split}
\right.
\end{equation} 
with
\begin{equation}
\label{d-wave}
d(t,x) = \left\{
\begin{split}
&\frac{\sigma(\sqrt{a(x)}z_t)}{\sqrt{a(x)}z_t},\quad &\text{ if } a(x)z_t\neq 0,\\
&C_1,\quad &\text{ if } a(x)z_t = 0.
\end{split}
\right. 
\end{equation}
This function comes from the function $t\mapsto \frac{\sigma(\sqrt{a(x)} z_t}{\sqrt{a(x)}z_t}$, extended at $0$, which is possible due to the differentiability of the function $\sigma$ at $0$.  

This system is in the form \eqref{ltv-system}, with $H=H_2$, $U=H_2$, $D(A) = D_2$, and the operators $A$ and $B$ defined as follows
\begin{equation}
\begin{split}
A : D(A) \subset H &\rightarrow H\\
\begin{bmatrix}
v_1 & v_2
\end{bmatrix}^\top &\mapsto \begin{bmatrix}
v_2 & v_1^{\prime\prime}
\end{bmatrix}^\top
\end{split}
\end{equation}
and $B=\begin{bmatrix}0 & \sqrt{a(x)} \end{bmatrix}^\top$. It is well known that $A$ generates a strongly continuous semigroup of contractions (see \cite{HarauxBook}). Now, let us check whether $d(t,\cdot)$ satisfies \eqref{bounds-D}, for all $t\geq 0$.

Since the initial conditions $(z_0,z_1)\in H_\infty(0,1)$ and are bounded by $R$ in the $H_\infty(0,1)$ norm, then invoking Theorem \ref{thm-wp} we have, for all $t\geq 0$
\begin{equation}
\sup_{x\in [0,1]} |z_t(t,x)|\leq 2R.
\end{equation} 
Moreover, since $\sqrt{a(x)}\leq \sqrt{a_\infty}$, for all $x\in [0,1]$, and because $d$ is a continuous function, there exist two positive constants $d_0$ and $d_1$, depending on $R$ and $a_\infty$ such that 
\begin{equation}
\begin{split}
d_0:=\min_{\xi\in [-2\sqrt{a_\infty} R,2\sqrt{a_{\infty}}R]} &\frac{\sigma(\xi)}{\xi} \leq d(t,x)\\
&\leq \max_{\xi\in [-2\sqrt{a_\infty} R,2\sqrt{a_{\infty}}R]} \frac{\sigma(\xi)}{\xi}:=d_1.
\end{split}
\end{equation}
Note moreover that the origin of the following system
\begin{equation}
\left\{
\begin{split}
&z_t = z_{xx} -d_0z_t\\
&z(t,0)=z(t,1)=0\\
&z(0,x)=z_0(x),\: z_t(0,x)=z_1(x).
\end{split}
\right.
\end{equation}
is exponentially stable in $H_2$ for any initial conditions $(z_0,z_1)\in H_2$ (see e.g. \cite{prieur2016wavecone}). The related operator of this system is $A-d_0BB^\star$, with domain $D(A_{d_0})=D_2$. Therefore, all the properties required in Theorem \ref{thm-ltv} are satisfied. Hence, there exist two positive constants $K:=K(R)$ and $\beta:=\beta(R)$ such that
\begin{equation}
\Vert (z,z_t)\Vert_{H_2} \leq K e^{-\beta t} \Vert (z_0,z_1)\Vert_{H_2},\quad \forall t\geq 0.
\end{equation}

\textbf{Second step: Semi-global exponential stability in $H_p(0,1)$}. From Theorem \ref{thm-wp}, we know that, for every initial conditions $(z_0,z_1)\in H_\infty(0,1)$ satisfying $\Vert (z_0,z_1)\Vert_{H_\infty(0,1)}\leq R$, and noticing that the trajectory of \eqref{wave-nonmonotone} can be expressed with the evolution family $W(t,0)$ associated to \eqref{wave-ltv}\footnote{Evolution families are extension of semigroups for time-variant linear infinite-dimensional systems.}, one has
\begin{equation}
\Vert W(t,0)(z_0,z_1)\Vert_{H_\infty(0,1)} \leq 2R, \quad \forall t\geq 0.
\end{equation}
Now, fix $t>0$. Note that $W(t,0)$ is a linear operator from $(L^2(0,1))^2$ (resp. $(L^\infty(0,1))^2$) to $(L^2(0,1))^2$ (resp. $(L^\infty(0,1))^2$), if it associates $(z_0^\prime,z_1)\in L^2(0,1)^2$ (resp. $(z_0^\prime,z_1)\in L^\infty(0,1)^2$)  to $(z_x,z_t)\in L^2(0,1)$ (resp. $(z_x,z_t)\in L^\infty(0,1)^2$). Hence, we can apply the so-called Riesz-Thorin theorem \cite[Theorem 1.1.1, Page 2]{bergh2012interpolation} and conclude that
\begin{equation}
\Vert W(t,0)\Vert_{H_p(0,1)}\leq 2\left(\frac{K}{2}\right)^{\frac{2}{p}}e^{-\frac{2\beta}{p}t}, \quad \forall t\geq 0,
\end{equation}
where $W(t,0)$ corresponds exactly to the trajectory of \eqref{wave-nonmonotone} with the initial condition $z_0,z_1\in H_\infty(0,1)$. In particular, for every $(z_0,z_1)\in H_\infty(0,1)$, one therefore has, for every $t\geq 0$
\begin{equation*}
\begin{split}
\Vert W(t,0)(z_0,z_1)\Vert_{H_p(0,1)} \leq &\Vert W(t,0)\Vert_{H_p(0,1)}\Vert (z_0,z_1)\Vert_{H_p(0,1)}\\
\leq & 2\left(\frac{K}{2}\right)^{\frac{2}{p}}e^{-\frac{2\beta}{p}t}\Vert (z_0,z_1)\Vert_{H_p(0,1)}.
\end{split}
\end{equation*}
This concludes the proof of Theorem \ref{thm-as}.\end{proof} 

\section{Conclusion}

\label{sec_conclusion}

In this paper, we have provided a well-posedness analysis of a one-dimensional wave equation subject to a nonlinear nonmonotone damping. Futhermore, a characterization of the asymptotic behavior of the latter system is given. It is proved with Lyapunov techniques.  This work paves the way to many others. For instance, a first open and natural question would be the case of multidimensional wave equations, for which there does not exist any proof of well-posedness in the functional setting introduced in this paper. 

\appendix
\section{Exponential convergence result for a linear time-variant system}
\label{sec_appendix}

This appendix is devoted to the statement and the proof of a theorem dealing with a time-variant linear infinite-dimensional system. Indeed, as it is illustrated in Section \ref{sec_as}, we can transform \eqref{wave-nonmonotone} as a time-variant linear infinite-dimensional system. This transformation is inspired by \cite{chitour1995continuity} and \cite{liu1996finite}, where such a technique is applied for finite-dimensional linear systems subject to a saturation. To define it, let us introduce $H$ (resp. $U$), a Hilbert space which is equipped with the norm $\Vert \cdot \Vert_H$ and the scalar product $\langle \cdot,\cdot\rangle_H$. The system under study in this section is the following:
\begin{equation}
\label{ltv-system}
\left\{
\begin{split}
&\frac{d}{dt} v = (A-d(t)BB^\star)v:=A_{d}(t)v\\
&v(\tau)=v_\tau,
\end{split}
\right.
\end{equation}
where $A:D(A)\subset H \rightarrow H$, with $D(A)$ the domain of the operator $A$ that we suppose densely defined in $H$, $B\in\mathcal{L}(U,H)$, $B^\star$ denotes the adjoint of $B$ and $d$ is a continuous function. We assume moreover that $A$ and its adjoint $A^\star$ are dissipative, and that there exist positive constants $d_0$ and $d_1$ such that $d$ satisfies, for all $t\geq 0$.
\begin{equation}
\label{bounds-D}
d_0 \leq d(t)\leq d_1.
\end{equation}
Due to these properties, the domain of $A_d(t)$ is equal to $D(A)$ and, therefore, does not depend on $t$. Moreover, always because of these properties, and since $A$ and $A^\star$ are dissipative, one can prove that, for all $t\geq 0$ and for all $v\in D(A)$
\begin{equation}
\langle A_d(t) v,v\rangle_H \leq 0,\: \langle A_d^\star(t) v,v\rangle_H\leq 0,
\end{equation}
which implies that each operator $A_d(t)$ generates a strongly continuous semigroup of contractions.  Hence, due to the discussion provided in \cite[Page 123]{chicone-latushkin1999}, there exists an evolution family $(W(\theta,\tau))_{\theta \geq \tau}$ solving \eqref{ltv-system}. In particular, this means that there exists a unique strong solution $v\in C(\mathbb{R}_+,D(A))$ to \eqref{ltv-system} for every $v_0\in D(A)$. 

We are now in position to state the following result.
\begin{theorem}
\label{thm-ltv}
Consider the system given by \eqref{ltv-system}. Supposing that the origin of the following system
\begin{equation}
\label{lti-system}
\left\{
\begin{split}
&\frac{d}{dt} v = (A-d_0 BB^\star) v:=A_{d_0}v,\\
&v(0)=v_0
\end{split}
\right.
\end{equation}
is globally exponentially stable, then, for any initial conditions $v_0\in H$, the origin of \eqref{ltv-system} converges exponentially to $0$ with $\tau=0$. 
\end{theorem}
\begin{remark}
The global asymptotic stability is not ensured for the linear time-variant system \eqref{ltv-system}. Indeed, for such systems, the global asymptotic stability has to be proven for every initial condition $v_\tau$ and for every $\tau\geq 0$. This is not an issue for us. Indeed, we use this result to prove the attractivity of \eqref{wave-nonmonotone}. The stability of \eqref{wave-nonmonotone} already holds thanks to Theorem \ref{thm-wp}.
\end{remark}

\begin{proof}
Since the origin of \eqref{lti-system} is globally exponentially stable, then, due to \cite{datko1970extending}, there exist a self-adjoint operator $P\in\mathcal{L}(H)$ and a positive constant $C$ such that the following inequality holds true
\begin{equation}
\langle PA_{d_0} v,v\rangle_H + \langle Pv,A_{d_0}v\rangle_H\leq -C\Vert v\Vert_H^2,\: \forall v\in D(A)
\end{equation}
Moreover, note that $A_{d_0}$ is also a dissipative operator, which means in particular that
\begin{equation}
\label{lyap-d0}
\langle A_{d_0} v,v\rangle_H + \langle v,A_{d_0}v\rangle_H \leq 0,\: \forall v\in D(A).
\end{equation}

Now, consider the following candidate Lyapunov functional for \eqref{ltv-system}:
\begin{equation}
V(v):=\langle Pv,v\rangle_H + M\Vert v\Vert_H^2,
\end{equation}
where $M$ is a positive constant which has to be defined. The time derivative of $V$ along the trajectories of \eqref{ltv-system} yields
\begin{equation}
\begin{split}
\frac{d}{dt} V(v)= &\langle (P+MI_H)A_{d(t)}v,v\rangle_H + \langle (P+MI_H)v,A_{d(t)} v\rangle_H\\
= & \langle PA_{d_0} v,v\rangle_H + \langle Pv,A_{d_0}v\rangle_H \\
 & + M(\langle A_{d_0} v, v\rangle_H + \langle v,A_{d_0}v\rangle_H)\\
 & -  \langle B(d(t) - d_0 I_H)B^\star v, (P+MI_H)v\rangle_H \\
 & -\langle (P+MI_H)v,B(d(t)-d_0I_H)B^\star v\rangle_H.\\
\leq & -C\Vert v\Vert_H^2 - d_0M\Vert B\Vert^2_{\mathcal{L}(H,U)}\Vert v\Vert^2_H\\
& + 2(d_1-d_0)\Vert B^\star\Vert_{\mathcal{L}(H,U)} \Vert P\Vert_{\mathcal{L}(H)} \Vert v\Vert_H^2 \\
&-2M \langle (d(t)-d_0)B^\star v,B^\star v\rangle_U,
\end{split}
\end{equation}
where we have used \eqref{bounds-D} and \eqref{lyap-d0} to get the inequality. Recalling that $d(t)-d_0\geq 0$, and setting 
$$
M = \frac{2(d_1-d_0)\Vert P\Vert_{\mathcal{L}(H)}}{d_0\Vert B\Vert_{\mathcal{L}(H,U)}},
$$
one has
\begin{equation}
\frac{d}{dt} V(v) \leq -C\Vert v\Vert_H^2,\quad \forall v\in D(A).
\end{equation}
Since $V$ satisfies the following inequalities
\begin{equation}
M\Vert v\Vert^2_H\leq V(z) \leq (\Vert P\Vert_{\mathcal{L}(H)} + M)\Vert v\Vert^2_H
\end{equation}
one can conclude that 
\begin{equation}
\Vert v\Vert_H^2 \leq \frac{\Vert P\Vert_{\mathcal{L}(H)} + M}{M}\exp\left(-\frac{C}{\Vert P\Vert_{\mathcal{L}(H)}+M}t\right) \Vert v_0\Vert_H^2,
\end{equation}
which ends the proof of Theorem \ref{thm-ltv}.
\end{proof}

\bibliography{bibsm}

\end{document}